\documentclass[11pt]{amsart}
\usepackage{enumerate}
\usepackage{fourier}
\usepackage{amssymb}
\usepackage{amsthm}
\usepackage{amsmath}
\usepackage{color}
\usepackage{enumerate}
\usepackage{fourier}
\usepackage{amssymb}
\usepackage{amsthm}
\usepackage{amsmath}
\usepackage{color}
\usepackage{amsmath}
\usepackage{amssymb}
\usepackage{amsfonts}
\usepackage{lipsum}
\usepackage{mathrsfs}

\usepackage{amsmath, cancel, amsthm, amscd, amsfonts, amssymb, graphicx, color}
\usepackage[bookmarksnumbered, colorlinks, plainpages]{hyperref}
\hypersetup{colorlinks=true,linkcolor=red, anchorcolor=green,
citecolor=cyan, urlcolor=red, filecolor=magenta, pdftoolbar=true}
\usepackage{ulem}

\newcommand{\B}{\mathcal B}

\newcommand{\C}{\mathbb{C}}

\newcommand{\h}{{\mathcal H}}

\newcommand{\lh}{{\mathcal B}({\mathcal H})}
\newcommand{\lx}{{\mathcal B}(X)}
\newcommand{\lxy}{{\mathcal B}(X,Y)}

\newcommand{\mn}{{\mathcal M}_n(\C)}

\newcommand{\uno}{{\bf\rm \textbf{1}}}

\DeclareMathSymbol{\subsetneqq}{\mathbin}{AMSb}{36}
\newtheorem{th1}{{\bfTheorem}}[section]
\newtheorem{thm}[th1]{{\bf Theorem}}

\newtheorem{lemma}[th1]{{\bf Lemma}}
\newtheorem{prop}[th1]{{\bf Proposition}}
\newtheorem{cor}[th1]{{\bf Corollary}}
\theoremstyle{remark}

\theoremstyle{definition}

\newtheorem{step}{Step}

\usepackage{color}

\usepackage{lipsum}

\newcommand{\Q}{\mathcal Q(X)}

\newcommand{\rank}{\operatorname{rank}}

\newcommand{\F}{\mathcal F}
\newcommand{\NN}{\mathcal N}

\def\B{{\mathcal B}(X)}
\def\Bh{{\mathcal B}({\mathcal H})}
\def\BB{{\mathcal B}(Y)}
\def\bxy{{\mathcal B}(Y,X)}
\def\bx*y{{\mathcal B}(Y,X^*)}
\def\byx*{{\mathcal B}(X^*,Y)}
\def\QX{{\mathcal Q}(X)}
\def\QY{{\mathcal Q}(Y)}
\def\FX{{\mathcal F}(X)}

\def\lxy{\mathcal{B}(X,Y)}

\def\lx{\mathcal{B}(X)}

\def\C{{\mathbb C}}
\def\r{{\rm r}}
\def\rank{{\rm rank}\,}

\title[Maps preserving quasi-nilpotent operators]
      {Surjective linear maps preserving quasi-nilpotent operators}

\date{}

\author{Abdellatif Bourhim }
\address{ Syracuse University, Department of Mathematics, 215 Carnegie Building, Syracuse, NY 13244, USA}
 \email{abourhim@syr.edu}
        
\author{Ajda Fo\v sner }
\address{Faculty of Management, University of Primorska, Koper, Slovenia}
 \email{ajda.fosner@fm-kp.si}

\author{Mostafa Mbekhta  }
\address{  Universit\'e de  Lille, 
          D\'epartement de Math\'ematiques,
        UMR-CNRS 8524, Labo. P. Painlev\'e \;  \quad
           59655 Villeneuve d'Ascq Cedex
        France.}

\email{mostafa.mbekhta@univ-lille.fr}


 \keywords{Bounded linear operator, Banach space, linear preserver, quasi-nilpotent operator, spectrum, spectral radius.\\
This work was partially supported in part by CDP C2EMPI and institutional
partners for the R-CDP-24-004-C2EMPI project and by Fincome program 470/2022 from CNRST- Morocco.}

\begin{document}

\maketitle{}



\begin{abstract}
Let $\Bh$ be the algebra of all bounded linear operators on an infinite-dimensional complex Hilbert space ${\mathcal H}$.
We show that a surjective linear map $\Phi:\Bh\to \Bh$ preserves quasi-nilpotent operators  if and only if there are a 
nonzero constant $c$ and an invertible operator $A\in\Bh$ such that either  $\Phi(T)=cATA^{-1}$ for all $T\in\Bh$ or 
$\Phi(T)=cAT^{tr}A^{-1}$ for all $T\in\Bh$. Here, $T^{tr}$ denotes the transpose of $T$ relative to a fixed but arbitrary 
orthonormal basis of ${\mathcal H}$. 
\end{abstract}

\section{Introduction}

Throughout this paper, $X$ and $Y$ denote two infinite-dimensional complex Banach spaces, and $\lxy$ denotes the space of all bounded linear maps from $X$ into $Y$. 
When $X=Y$, we simply write $\lx$ instead of ${\mathcal B}(X,X)$ and denote its identity operator by  $\uno\in \B$. Let  $\NN(X)$  be  the set of all nilpotent operators of $\lx$, and ${\mathcal B}_0(X)$ be its linear span. As usual, the dual space of $X$ will be denoted by $X^*$ and the Banach space adjoint of any operator $T\in \B$ will be denoted by $T^*$. In \cite{BresarSemrl}, Bre\v sar and \v Semrl proved that every surjective linear map on $\B$ preserving the spectral radius is a Jordan automorphism multiplied by a complex constant of modulus one. Their work was motivated by Jafarian and Sourour's result on spectrum-preservers \cite{JS} and, more broadly, by Kaplansky's problem \cite{Kaplanski} on possible generalizations of the well known Gleason-Kahane-\. Zelazko theorem; see for instance \cite{Aupetit1,Aupetit2,Aupetit3}. The literature surrounding the aforementioned results is too extensive to review here, and we refer the interested reader to the survey article \cite{survey} and the references therein. We nevertheless single out  a \v Semrl's  result, as it is closely related to the problem that will be addressed in this paper.
\begin{thm}[\cite{Semrl1}]\label{nil-Semrl}
If $\Phi:{\mathcal B}_0(X)\longrightarrow{\mathcal B}_0(Y)$ is a surjective linear map preserving nilpotent operators, then  
there exist a nonzero scalar $c\in\C$, and either a bijective map $A\in\bxy$ such that 
\begin{equation}
\label{sform1}
\Phi(T)=cATA^{-1},\quad \left(T\in{\mathcal B}_0(X)\right),
\end{equation}
or a bijective map $B\in\byx*$ such that 
\begin{equation}\label{sform2}
\Phi(T)=cBT^*B^{-1},\quad \left(T\in{\mathcal B}_0(X)\right).
\end{equation}
In the last case, both spaces $X$ and $Y$ must be reflexive.
\end{thm}

Closely related to this result is a longstanding problem concerning the description of surjective linear maps on $\lx$ preserving the quasi-nilpotent operators. In view of the latter result, one may wonder whether analogous conclusions can be established for such a map. Although several attempts have been made to resolve this problem, none has succeeded in establishing the desired conclusion in its full generality; see for instance \cite{costara}. The purpose of this paper is to settle this problem by demonstrating that an analogous conclusion indeed holds. 

Our investigation of this problem dates back to 2011, when we developed a proof that remained incomplete due to a single missing step, namely, establishing that the rank-one nilpotent operators are preserved in both directions by any surjective linear map on $\lx$ preserving the quasinilpotency. Despite our efforts at the time, this identity remained elusive. More recently, with the assistance of artificial intelligence (ChatGPT), we were able to establish this crucial identity, thereby bridging the gap in our original argument and obtaining a complete proof of the result.

\section{Main result}
For any operator $T\in\lx$, let $\sigma(T)$ be its spectrum and by $\r(T)$ be its spectral radius of $T$.  Recall that an operator $T\in \B$ is said to be quasi-nilpotent if $\r(T)=0$ (or equivalently, $\sigma(T)=\{0\}$). Obviously, every nilpotent operator is quasi-nilpotent; however, a quasi-nilpotent operator need not be nilpotent. Throughout, we denote by $\Q$  the set of all quasi-nilpotent operators    of $\lx$, and say that a linear map $\Phi : \B \to \BB$ preserves quasi-nilpotent operators if, for every $T\in \B$, we have  
$$T\in \QX \iff \Phi(T)\in \QY.$$ In other words, $\Phi$ preserves quasi-nilpotent operators if $$\r(T)=0 \iff \r(\Phi(T))=0$$ for every $T\in \B$. 

The following result is the main result of this paper and it describes all surjective linear maps on operator algebras preserves quasi-nilpotent operators.  It shows that such a map takes one of the above forms (\ref{sform1}) or (\ref{sform2})
when restricted to ${\mathcal B}_0(X)$.

\begin{thm}\label{main}
If $\Phi:\B\to \BB$ is a surjective linear map preserving quasi-nilpotent operators, then $\Phi$ maps ${\mathcal B}_0(X)$ onto ${\mathcal B}_0(Y)$ and 
there exist a nonzero scalar $c\in\C$, and either a bijective map $A\in\bxy$ such that 
\begin{equation}
\label{form1}
\Phi(T)=cATA^{-1},\quad \left(T\in{\mathcal B}_0(X)\right),
\end{equation}
or a bijective map $B\in\byx*$ such that 
\begin{equation}\label{form2}
\Phi(T)=cBT^*B^{-1},\quad \left(T\in{\mathcal B}_0(X)\right).
\end{equation}
In the last case, both spaces $X$ and $Y$ must be reflexive.
\end{thm}

Note that if $X=Y$ is a finite-dimensional Banach space and $\dim (X) = n$, the algebra $\B$ can be canonically identified with the matrix algebra $\mn$ of all $n\times n$-complex matrices. By the Cayley-Hamilton theorem, every quasi-nilpotent matrix in $\mn$ is nilpotent, and the characterization of all surjective linear maps on $\mn$ preserving nilpotent matrices  is already known and is due to Botta et al.; see \cite{Botta}. This justifies our assumption throughout the present paper that the Banach spaces $X$ and $Y$  are assumed to be infinite-dimensional.

We close this section by stating an immediate consequence of Theorem \ref{main}.  To do so, let $\h$ be an infinite-dimensional complex Hilbert space, and note that, since every operator in $\Bh$ is a finite sum of square-zero operators \cite{CPT, Fillmore, WWu}, we see that ${\mathcal B}_0(\h)=\Bh$.
\begin{cor}
A surjective linear map $\Phi:\Bh\to\BB$ preserves quasi-nilpotent operators if and only if there exist  a bijective operator
$A\in{\mathcal B}({\mathcal H},Y)$ and a nonzero scalar $c$ such that either $\Phi(T)=cATA^{-1}$ for all $T\in\Bh$, or $\Phi(T)=cAT^{\text{tr}}A^{-1}$ for all $T\in\Bh$.
Here, $T^{tr}$ denotes the transpose of any operator $T\in\lh$ with respect to an arbitrary fixed basis of ${\mathcal H}$.  
\end{cor}

\section{Preliminaries and auxiliary  results}
In this section we gather together some results and notation needed for the proof of Theorem \ref{main}.
For any  nonzero vector $x\in X$ and nonzero linear functional $f\in X^*$, let  $$(x\otimes f)(z):=f(z)x, \quad (z\in X).$$ 
Note that $x\otimes f$ is an operator of rank one and that any rank-one operator in $\lx$ is of this form. Moreover, $x\otimes f$ is nilpotent if and only if $f(x)=0$. Furthermore, every finite-rank operator $T\in\lx$ can be written as a finite sum of  rank-one operators in $\lx$; i.e., 
$T=\sum_{k=1}^nx_k\otimes f_k$ for some $x_k\in X$ and $f_k\in X^*$, $k=1,2,\ldots , n$. Let ${\mathcal F}_1(X)$, $\NN_1(X)$ and  $\FX$ stand for the set of all rank-one operators, for the set of all nilpotent rank-one operators, and for the set of all finite rank operators in $\B$, respectively. 
\medskip

We first record an elementary yet useful lemma that identifies those points of the spectrum of a rank-one perturbation of an operator $T\in\lx$ that lie outside the spectrum of $T$.

\begin{lemma}
\label{L:sourour1}
Let $T\in\lx$ be an operator and let $\lambda\in\C\setminus\sigma(T)$. If $x\in X$ and  $f\in X^*$, then 
  $\lambda\in\sigma(T+x\otimes f)$ if and only if $f((\lambda-T)^{-1}x)=1$.
\end{lemma}
\begin{proof}
See, for example, \cite[Lemma 4]{JS}.
\end{proof}

The second lemma is well known and its proof is quoted from \cite{BresarSemrl} for the sake of completeness of the reader. 

\begin{lemma}[Bre\v sar--\v Semrl \cite{BresarSemrl}]\label{zemanek}
An operator $T_0\in\lx$ is identically zero if and only if $\r(T_0+Q)=0$ for all  $Q\in\Q$.
\end{lemma}
\begin{proof}
Assume that $\r(T_0+Q)=0$ for all  $Q\in\Q$ and suppose for the sake of contradiction that $T_0\ne 0$. Then  there is a vector $x$ in $X$ such that $T_0x\not=0$.  Since $\r(T_0)=0$, the vectors $T_0x$ and $x$ are linearly independent. Hence, there is a linear functional $f\in X^*$ such that
$f(x)-=f(T_0x)=1$.  Set $Q:=(-T_0x+x)\otimes f$ and note that $Q^2=0$.  Since $(T_0+Q)x=x$, we get  $\r(T_0+Q)\geq1$, which is a contradiction. This complete the proof.
\end{proof}

We say that a linear map $\Phi:{\mathcal B}(X)\to {\mathcal B}(Y)$ preserves nilpotent rank-one operators in both directions if for every $T\in \B$ we have
$$T\in \NN_1(X) \iff \Phi(T)\in \NN_1(Y).$$ The next well-known result characterizes bijective linear maps 
between ${\mathcal F}(X)$  and  ${\mathcal F}(Y)$ which preserve rank-one operators (see, for example, \cite{OS}).

\begin{thm}
\label{T2}
If $\Phi:{\mathcal B}(X)\to {\mathcal B}(Y)$ is a bijective linear map preserving rank-one nilpotent operators in both directions,  
then there exist a nonzero scalar $c\in\C$ and either a bijective map $A\in\bxy$ such that
\begin{equation}
\label{1}
\Phi(T)=cATA^{-1},~(T\in\NN_1(X)),
\end{equation}
or a bijective map $B\in\byx*$ such that 
\begin{equation}
\label{2}
\Phi(T)=cBT^*B^{-1},~(T\in\NN_1(X)).
\end{equation}
The latter case may occur only if both spaces $X$ and $Y$ are reflexive.
\end{thm}
\begin{proof}
See for instance \cite{OS} or \cite{Du}.
\end{proof}

For each operator in $\lx$, we now introduce a family of sets that will play a crucial role in what follows. Indeed, for any operator $T\in\lx$ and $Q\in\Q$, let 
$$\Omega_T(Q):=\left\{\lambda\in\C\setminus\{0\}:r(Q+\lambda T)=0\right\},$$
and denote the cardinality of any subset $\Lambda$ of $\C$ by $\left|\Lambda\right|$. Our purpose of the rest of this section is to make use of this set and establish a couple of interesting results in their own right that are need for the proof of the main result. We start with the following lemma that tells us that if $N\in\NN_1(X)$ then  the cardinal of each $\Omega_N(Q)$ is no more that $1$ whenever $Q\in\Q\backslash\{0\}$ and  $\r(Q+N)>0$.

\begin{lemma}\label{lem:rankone}
If $N\in\NN_1(X)$, then $$\left|\Omega_N(Q)\right|\leq 1$$
for any nonzero $Q\in\Q$ for which $\r(Q+N)>0$.
\end{lemma}
\begin{proof}Let $N\in\NN_1(X)$ and $x\in X$ be a nonzero vector and $f\in X^*$ be a nonzero linear functional such that $f(x)=0$ and $N=x\otimes f$. Let  $Q\in\Q$ be a quasi-nilpotent operator such that $\r(Q+N)>0$, and let
$$
F(z):=\sum_{n=0}^{\infty} f(Q^nx)z^{n+1}.
$$
Since $0=\r(Q)=\lim_{n\to\infty}\|Q^n\|^{1/n}$, this  series converges locally uniformly on $\C$ and then $F$ is an entire function. Note that, since  $0<\r(Q+N)=\r(Q+x\otimes f)$, there exists a nonzero $\lambda_0\in\sigma(Q+x\otimes f)$. Then Lemma \ref{L:sourour1} tells us that 
$$1=f((\lambda_0\uno-Q)^{-1}x)=\sum_{n\geq0}\frac{f(Q^nx)}{\lambda_0^{n+1}}=F\left(\lambda_0^{-1}\right).$$
This and the fact that $F(0)=0$ entail  that $F$ is a nonconstant entire function. Consequently, by Great Picard's theorem, $F$ assumes all complex values except possibly, say $z_0\in\C$. For any nonzero $z\in\C$ for which $zz_0\not=1$, there is a nonzero $\lambda\in\C$ such that $\frac{1}{z}=F(\lambda)=f\left((\lambda^{-1}\uno-Q)^{-1}x\right),$
and thus Lemma \ref{L:sourour1} implies that 
$\frac1\lambda\in\sigma\left(Q+zN\right)$. Clearly, $\r(Q+z N)>0$ for all nonzero $z\in\C$ except perhaps for $z=\frac1{z_0}$. Therefore, $\Omega_N(Q)\subseteq\Bigl\{\frac1{z_0}\Bigr\}$ and 
$$\left|\Omega_N(Q)\right|\leq 1$$
for any nonzero $Q\in\Q$ for which $\r(Q+N)>0$.
This completes the proof.
\end{proof}

In the remainder of this section, we investigate the converse of the preceding lemma. Our next result asserts that every nonzero quasi-nilpotent operator  satisfying the conclusion of Lemma \ref{lem:rankone}  is nilpotent and its index of nilpotency at most $5$. 

\begin{prop}\label{prop:five}Let  $R\in\Q$ be a nonzero quasi-nilpotent operator. If 
$$\left|\Omega_R(Q)\right|\leq 1$$
for any nonzero operator $Q\in\Q$ of rank at most $3$ and of index of nilpotency at most $4$ for which $\r(Q+R)>0$, then $R^5=0$.
\end{prop}
\begin{proof}
Suppose for the sake of contradiction that $R^5\ne 0$, and let $u\in X$ such that $R^5u\ne 0$.
Since $R$ is quasinilpotent,  the vectors $u,~Ru,~R^2u,~R^3u,~R^4u,~R^5u$
are linearly independent. Hence, there exists $f\in X^*$
such that
\begin{equation}\label{choice}
f(u)=f(Ru)=0,\qquad
f(R^ju)=1,\quad j=2,3,4,5.
\end{equation}
Set $x_0:=u-Ru$ and note that 
\begin{equation}\label{frx}
f(Rx_0)=-1, \qquad f(x_0)=f(R^jx_0)=0,\qquad j=2,3,4.
\end{equation}
Let 
\begin{equation}\label{defQ}
Q:=(x_0\otimes f) R^{2}+(Rx_0\otimes f)R +(R^2x_0\otimes f),
\end{equation}
and note that $Q$ is of rank at most $3$. In view of \eqref{frx}, we see that 
\begin{equation}\label{frxbis}
Qx_0=-Rx_0,\qquad QRx_0=-R^2x_0,\qquad QR^2x_0=0.
\end{equation}
Consequently,$$
\begin{aligned}
Q^2&=Q\left((x_0\otimes f) R^{2}+(Rx_0\otimes f)R +(R^2x_0\otimes f)\right)\\
&=(Qx_0\otimes f) R^{2}+(QRx_0\otimes f)R +(QR^2x_0\otimes f)\\
&=-(Rx_0\otimes f) R^{2}-(R^2x_0\otimes f)R,\\[1ex]
Q^3&=-Q\left((Rx_0\otimes f) R^{2}+(R^2x_0\otimes f)R\right)\\
&=(R^2x_0\otimes f) R^{2},\text{ and }\\
Q^4&=Q\left((R^2x_0\otimes f) R^{2}\right)=0.
\end{aligned}
$$
Clearly, $Q\in\Q$ is an operator of rank at most $3$ and its index of nilpotency is at most $4$. On the other hand, in view of \eqref{choice}, we obtain
$$Qu=f(R^2u)x_0+f(Ru)Rx_0+f(u)R^2x_0=x_0=u-Ru,$$
Thus, $(R+Q)u=u$ and  
\begin{equation}\label{eq:positive-radius}
\r(R+Q)\geq 1.
\end{equation}

Now, fix a nonzero $s\in\C$, let 
$$S_s:=\uno+(1-s)(x_0\otimes f)R+(1-s^2)Rx_0\otimes f.$$
A direct computation shows that $S_s$ is invertible and its inverse is $S_{s^{-1}}.$ In view of \eqref{defQ} and \eqref{frxbis}, we also have 
\begin{eqnarray*}
QS_s&=&Q+(1-s)(Qx_0\otimes f)R+(1-s^2)(QRx_0)\otimes f\\
&=&Q-(1-s)(Rx_0\otimes f)R-(1-s^2)R^2x_0\otimes f\\
&=&(x_0\otimes f) R^{2}+s(Rx_0\otimes f)R+s^2R^2x_0\otimes f.
\end{eqnarray*}
Moreover, we have 
$$RS_s=R+(1-s)(Rx_0\otimes f)R+(1-s^2)R^2x_0\otimes f,$$
and
$$S_sR=R+(1-s)(x_0\otimes f)R^2+(1-s^2)(Rx_0\otimes f)R.$$
Then 
\begin{eqnarray*}
RS_s-S_sR&=&(s^2-s)(Rx_0\otimes f)R+(1-s^2)R^2x_0\otimes f-(1-s)(x_0\otimes f)R^2\\
&=&(s-1)\left(s(Rx_0\otimes f)R-(1+s)R^2x_0\otimes f+(x_0\otimes f)R^2\right)\\
&=&(s-1)\left(QS_s-(1+s+s^2)R^2x_0\otimes f\right)\\
&=&(s-1)QS_s+(1-s^3)R^2x_0\otimes f.
\end{eqnarray*}
Whence 
\begin{equation}\label{eq:intertwining}
\Bigl(R+(1-s)Q\Bigr)S_s=S_sR+(1-s^3)R^2x_0\otimes f
\end{equation}
for all nonzero $s\in\C$.

Let $\omega$ be a primitive cube root of unity.
Taking $s=\omega$ and $s=\omega^2$, respectively,
in \eqref{eq:intertwining}, we obtain
\[
(R+(1-\omega)Q)S_\omega=S_\omega R
\]
and
\[
(R+(1-\omega^2)Q)S_{\omega^2}
=S_{\omega^2}R.
\]
Since $S_\omega$ and $S_{\omega^2}$ and  $R$ is quasinilpotent, we conclude that
$$\r(R+(1-\omega)Q)=\r(R+(1-\omega^2)Q)=0.$$
Since $1-\omega$ and $1-\omega^2$ are distinct
nonzero scalars and their reciprocals are in $\Omega_R(Q)$, we see that 
$$\left|\Omega_R(Q)\right|\geq 2.$$ 
This contradiction proves that $R^5=0$, and completes the proof.
\end{proof}

The following result is inspired by the proof
of Proposition~2.1 in \cite{Semrl1},
with suitable modifications to accommodate two
distinct exceptional nonzero scalars.

\begin{prop}\label{prop:rank}
Let $R\in\lx$ be  a nonzero nilpotent operator. If 
$$\left|\Omega_R(Q)\right|\leq 1$$
for any nonzero operator $Q\in\Q$ of rank at most $3$ and of index of nilpotency at most $4$ for which $\r(Q+R)>0$, then $R$ is of rank at most $2$. Consequently, $R^3=0$.
\end{prop}

\begin{proof}
Suppose, toward a contradiction, that the rank of $R$ is at least 3 and let $n$ be a positive integer such that $R^n=0$.
Choose $x_1,x_2,x_3\in X$ such that the vectors
$Rx_1,Rx_2,Rx_3$ are linearly independent. Let 
$$
X_0:=\operatorname{span}
\Bigl\{R^jx_i:0\leq j<m,\ 1\leq i\leq 3\Bigr\}.
$$
Clearly, $X_0$ is a finite-dimensional invariant
subspace of $R$ and $$k:=\rank(R|_{X_0})\geq 3.$$
Let $d:=\dim(X_0)$, and note that, since $R|_{X_0}$ is nilpotent,
the Jordan canonical form theorem ensures the
existence of an ordered basis
$\{e_1,\ldots,e_d\}$ of $X_0$ such that
\[
Re_i=a_i e_{i+1},\qquad i=1,\ldots,d,
\]
where $e_{d+1}=e_1$ and
\[
a_i\in\{0,1\},\qquad
\sum_{i=1}^d a_i=k.
\]
By the Hahn--Banach theorem, there exist
linear functionals $f_1,\ldots,f_d\in X^*$ such that
\[
f_i(e_j)=\delta_{ij},
\qquad 1\leq i,j\leq d.
\]
Since $k\geq 3$, we can choose two distinct indices
$p,q$ such that $a_p=a_q=1$. Following the the arguments of 
\cite[Proof of Proposition~2.1]{Semrl1}, define the
finite-rank operator $Q\in\lx$ by
\begin{equation}\label{eq:rank-test}
Q:=-2e_{p+1}\otimes f_p
  -3e_{q+1}\otimes f_q
  +\sum_{\{i:a_i=0\}}e_{i+1}\otimes f_i.
\end{equation}
Clearly, $Q(X)\subseteq X_0$. Moreover,
\[
Qe_i=b_i e_{i+1},
\qquad i=1,\ldots,d,
\]
where
\[
b_i=
\begin{cases}
-2&\text{if }i=p,\\
-3&\text{if }i=q,\\
1&\text{if }a_i=0,\\
0&\text{otherwise}.
\end{cases}
\]
Now, for every $\lambda\in\C$, set
\[
M_\lambda:=(Q+\lambda R)|_{X_0}.
\]
Then
\begin{eqnarray*}
M_\lambda e_i
&=&Qe_i+\lambda Re_i\\
&=&(b_i+\lambda a_i)e_{i+1}.
\end{eqnarray*}
Applying $M_\lambda$ successively $d$ times,
we obtain
\begin{eqnarray*}
M_\lambda^d e_i
&=&\left(\prod_{j=1}^d
(b_j+\lambda a_j)\right)e_i\\
&=&\lambda^{k-2}(\lambda-2)(\lambda-3)e_i.
\end{eqnarray*}
Since $\{e_1,\ldots,e_d\}$ is a basis of $X_0$,
it follows that
\begin{equation}\label{eq:semrl-cycle}
M_\lambda^d
=\lambda^{k-2}(\lambda-2)(\lambda-3)\uno_{X_0}.
\end{equation}
Since $k\geq 3$, taking $\lambda=0$ in
\eqref{eq:semrl-cycle} yields
\[
(Q|_{X_0})^d=0.
\]
Together with the fact that $Q(X)\subseteq X_0$,
this implies that $Q^{d+1}=0$. Thus $Q$ is a
finite-rank nilpotent operator.

Taking $\lambda=2$ and $\lambda=3$ in
\eqref{eq:semrl-cycle}, we obtain
\[
M_2^d=M_3^d=0.
\]
On the other hand, since $Q(X)\subseteq X_0$ and
$R(X_0)\subseteq X_0$, a straightforward induction
shows that
\[
(Q+\lambda R)^j x-\lambda^j R^j x\in X_0
\]
for every $x\in X$ and every positive integer $j$.
Since $R^n=0$, we obtain
\[
(Q+\lambda R)^n(X)\subseteq X_0.
\]
Consequently, for $\lambda=2,3$, we have
\begin{eqnarray*}
(Q+\lambda R)^{n+d}
&=&(Q+\lambda R)^d(Q+\lambda R)^n\\
&=&0.
\end{eqnarray*}
Therefore,
\begin{equation}\label{eq:semrl-zero}
\r(Q+2R)=\r(Q+3R)=0,
\end{equation}
and 
$$\left|\Omega_R(Q)\right|\geq 2.$$ 
To reach a contradiction, we only need to show that $\r(Q+R)>0$. To do so,  take $\lambda=1$ in
\eqref{eq:semrl-cycle}, and note that 
\[
M_1^d=2\uno_{X_0}.
\]
Consequently,
\begin{eqnarray*}
\r(M_1)^d
&=&\r(M_1^d)\\
&=&\r(2\uno_{X_0})\\
&=&2,
\end{eqnarray*}
and hence
\[
\r(M_1)=2^{1/d}>0.
\]
Since $X_0$ is a finite-dimensional invariant
subspace of $Q+R$, every eigenvalue of $M_1$
is also an eigenvalue of $Q+R$. Therefore,
\begin{equation}\label{eq:semrl-positive}
\r(Q+R)\geq r(M_1)=2^{1/d}>0.
\end{equation}
Therefore, the desired contradiction is established, implying that the rank of $R$ is at most $2$, which completes the proof.
\end{proof}

The last lemma characterizes the nilpotent rank-one operators  of $\lx$ and distinguishes them from the finite-rank nilpotent operators with higher ranks.
\begin{lemma}\label{strong}
A nonzero finite-rank nilpotent operator $R\in\lx$ has rank
one if and only if $\Omega_R(N)=\emptyset$ for all  finite-rank nilpotent
operators $N\in\NN(X)$ with $r(N+R)>0$.
\end{lemma}
\begin{proof} If $R$ has rank one and $N\in\NN(X)$ is a finite-rank nilpotent operator for which  $\r(N+R)>0$, then \cite[Proposition 2.1]{Semrl1} tells us that 
$N+\lambda R\not\in\NN(X)$ for all nonzero $\lambda\in\C$. Since each $N+\lambda R$ is a finite-rank operator, we see that 
$$\r(N+\lambda R)>0$$
for all nonzero $\lambda\in\C$, and thus  $\Omega_R(N)=\emptyset$. 

Conversely, suppose that the rank of $R$ is at least $2$.
By the construction used  in \cite[Proof of
Proposition~2.1]{Semrl1}, there exists
a finite-rank nilpotent operator $N\in\NN(X)$ such that
$N+R$ is not nilpotent, whereas $N+2R$ is
nilpotent. Since both $N+R$ and $N+2R$ are finite-rank
operators, it follows that
\[
r(N+R)>0
\quad\text{and}\quad
r(N+2R)=0,
\]
and thus $2\in\Omega_R(N)$. This contradicts the hypothesis of this lemma and completes the proof.
\end{proof}

\section{Proof of the main theorem}

Assume that $\Phi:\B\to \BB$ is a surjective linear map preserving quasi-nilpotent operators, and let us  divide the proof into several steps.
\begin{step}
 $\Phi$ is injective.
\end{step}

Assume that $\Phi(T_0)=0$ for some $T_0\in\B$. For every quasi-nilpotent operator $Q\in \B$ we have
$$\r(Q)=0\iff\r(\Phi(T_0)+\Phi(Q))=0\iff\r(T_0+Q)=0.$$ Thus, Lemma \ref{zemanek} tells us that $T_0=0$ and the map $\Phi$ is injective. Since $\Phi$ is assumed to be surjective, the map $\Phi$ is, in fact, bijective.

\begin{step}
$\Phi\left(\F_0(X)\right)=\F_0(Y)$. 
\end{step}

Assume that $N\in\NN_1(X)$ and note that Lemma \ref{lem:rankone} implies that $\left|\Omega_N(Q)\right|\leq 1$
for any nonzero $Q\in\Q$ for which $\r(Q+N)>0$. Thus, 
$$\left|\Omega_{\Phi(N)}(\Phi(Q))\right|\leq 1$$
for any nonzero $Q\in\Q$ for which $\r(Q+N)>0$. Since $\Phi$ is bijective and preserves the quasi-nilpotency in both direction,  Proposition \ref{prop:five} and Proposition \ref{prop:rank} tell us that $R$ is a nilpotent operator of rank at most two and $\Phi(N)\in \F_0(Y)$. Consequently,  $\Phi\left(\NN_1(X)\right)\subseteq\F_0(Y).$
Since $\NN_1(X)$ spans $\F_0(X)$ and $\Phi$ is linear, we see that 
$$\Phi\left(\F_0(X)\right)\subseteq\F_0(Y).$$ 
Since $\Phi^{-1}$ preserves the quasi-nilpotency in both directions, what has been shown previously applied to $\Phi^{-1}$ yields the reverse inclusion and thus $\Phi\left(\F_0(X)\right)=\F_0(Y)$.

\begin{step}
There exist a nonzero scalar $c\in\C$, and either a bijective map $A\in\bxy$ such that 
\begin{equation}
\label{sssform1}
\Phi(T)=cATA^{-1},\quad \left(T\in\NN_1(X)\right),
\end{equation}
or a bijective map $B\in\byx*$ such that 
\begin{equation}\label{sssform2}
\Phi(T)=cBT^*B^{-1},\quad \left(T\in\NN_1(X)\right).
\end{equation}
\end{step}

Let $$\NN_{f}(X):=\F(X)\cap\NN(X)$$
be the set of all finite-rank nilpotents and define $\NN_f(Y)$ in a similar way. Note that 
$$\NN_{f}(X)=\F_0(X)\cap\Q.$$
In view of the previous step, we have $\Phi\left(\F_0(X)\right)=\F_0(Y)$. By the preserving property of $\Phi$, we also have 
$\Phi\left(\Q\right)={\mathcal Q}(Y).$
As $\Phi$ is bijective, we obtain 
\begin{eqnarray*}
\Phi(\NN_f(X))&=&\Phi\left(\F_0(X)\cap\Q\right)\\
&=&\Phi\left(\F_0(X)\right)\cap\Phi\left(\Q\right)\\
&=&\F_0(Y)\cap{\mathcal Q}(Y)\\
&=&\NN_f(Y).
\end{eqnarray*}
This and Lemma \ref{strong} entail that $\Phi(\NN_1(X))=\NN_1(Y)$, and consequently, $\Phi$ maps ${\mathcal B}_0(X)$ onto ${\mathcal B}_0(Y)$. Therefore, by Theorem \ref{T2}, the map $\Phi$ takes either the form \eqref{sssform1} or  \eqref{sssform2}

\begin{step}
Without loss of generality, we may and will assume that $$\Phi(T)=T$$ for all $T\in\NN_1(X)$.
\end{step}
If $\Phi$ takes the form \eqref{sssform2}, then $Y$ must be reflexive and thus  the natural embedding $K:Y\to Y^{**}$ is bijective. Replacing $\Phi$ by $T\mapsto \frac{1}{c}A^{-1}\Phi(T)A$ when $\Phi$ takes the form \eqref{sssform1} and  by $T\mapsto \frac{1}{c}K^{-1}B^*\Phi(T)^*(B^{-1})^*K$ if $\Phi$ takes the form \eqref{sssform2}. Thus, there is no loss of generality if $$\Phi(T)=T$$for all $T\in\NN_1(X)$.
Our aim is then  to prove that $\Phi(T)=T$ for all $T\in{\mathcal B}_0(X)$. Thus, it suffices to show that $\Phi(N)=N$ for all nilpotent operators $N\in\B$ since $\Phi$ is linear.

\begin{step}
If $Q\in\B$ is a quasi-nilpotent operator and $Qx=0$ for some $x\in X$, then $\Phi(Q)x=0$.
\end{step}

Assume that $Qx=0$ for some $x\in X$. It suffices to show that $x$ and $\Phi(Q)x$ are linearly dependent since $\Phi(Q)$ is a quasi-nilpotent operator.
Assume on the contrary that they are linearly independent. Then $x$ and $u:=(I - \Phi (Q))^{-1}x$ are linearly independent, as well. Let us pick a linear 
functional $f\in X^*$ such that $f(x) = 0$ and $ f(u) =1$. Then 
\begin{eqnarray*}
\big(\Phi (Q) + x\otimes f\big) u &=&\big((\Phi (Q) - I) + I + x\otimes f\big)u\\
&=&(\Phi (Q) - I)u+u+f(u)x\\
&=&-x+u+x=u.
\end{eqnarray*}
This implies that $\r(\Phi(Q) + x\otimes f)>0$. On the other hand, since $Qx =0$ and $f(x) =0$, we have 
$$f((\lambda - Q)^{-1}x)=\sum_{n\geq0}\frac{f(Q^nx)}{\lambda^{n+1}}=0$$ for all $\lambda\not=0$.
This shows that  $\lambda - Q - x\otimes f = (\lambda - Q)(I - (\lambda - Q)^{-1}x \otimes f)$ is invertible for all $\lambda\not=0$. Thus, $\r(Q+x\otimes f)=0$ 
and, consequently, $\r(\Phi (Q + x\otimes f)) = \r(\Phi (Q) +x\otimes f) =0$. This contradiction shows that $x$ and $\Phi(Q)x$ are linearly dependent and, thus, $\Phi(Q)x=0$, as desired.

\begin{step}
$\Phi(N)=N$ for all nilpotent operators $N\in\B$.
\end{step}

Let $N$ be a nilpotent operator and $x\in X$. We would like to show that $\Phi(N)x=Nx$. By the previous step, we may assume that $Nx\not=0$. Let $k$ be 
the largest integer for which $N^kx\not=0$. Note that $x,Nx,\dots,N^kx$ are linearly  independent. Let us choose any $f\in X^*$ such that $f(x)=1$ and 
$f(N^nx)=0$ for $1\leq n\leq k$. Then $$f((\lambda - N)^{-1}Nx)=\sum_{n\geq0}\frac{f(N^{n+1}x)}{\lambda^{n+1}}=0$$ for all $\lambda\not=0$. As above, 
this shows that $N-Nx\otimes f$ is quasi-nilpotent. Since $(N-Nx\otimes f)x=0$, the previous step implies that $\Phi(N-Nx\otimes f)x=0$ and, thus, $\Phi(N)x=Nx$.
In other words, $\Phi(N)=N$ for all nilpotent operators $N\in\B$. This completes the proof.

\end{document}